\documentclass[10pt, reqno]{amsart}
\usepackage{amsmath, amsthm, amscd, amsfonts, amssymb, graphicx, color}
\usepackage[bookmarksnumbered, colorlinks, plainpages]{hyperref}
\hypersetup{colorlinks=true,linkcolor=red, anchorcolor=green, citecolor=cyan, urlcolor=red, filecolor=magenta, pdftoolbar=true}
\usepackage{mathrsfs}

\newtheorem{theorem}{Theorem}[section]
\newtheorem{lemma}[theorem]{Lemma}
\newtheorem{proposition}[theorem]{Proposition}

\theoremstyle{definition}

\theoremstyle{remark}
\newtheorem{remark}[theorem]{Remark}
\numberwithin{equation}{section}

\begin{document}
\setcounter{page}{1}

\title[Hilbert 12th problem]{Measured foliations and Hilbert 12th problem}

\author[Nikolaev]
{Igor V. Nikolaev$^1$}

\address{$^{1}$ Department of Mathematics and Computer Science, St.~John's University, 8000 Utopia Parkway,  
New York,  NY 11439, United States.}
\email{\textcolor[rgb]{0.00,0.00,0.84}{igor.v.nikolaev@gmail.com}}





\begin{abstract}
Yu.~I.~Manin conjectured  that the maximal abelian extensions of 
the real quadratic number fields are generated by the  pseudo-lattices
with real multiplication. We prove this   conjecture   using  theory of 
measured  foliations on the modular curves.
\end{abstract}

\maketitle

\section{Introduction}
The Hilbert 12th problem consists in a generalization  of the Kronecker-Weber theorem for  abelian extensions of the rational numbers to any base number field.  It therefore asks for analogues of the roots of unity, and of their appearance as particular values of a special function, in this case the exponential function;  the requirement is that such numbers should generate a whole family of further number fields that are analogues of the cyclotomic fields and their subfields.  The classical theory of complex multiplication solves this problem for the case of any imaginary quadratic field, by using modular functions and elliptic functions with respect to a particular period lattice related to the field in question. Shimura extended this to general CM fields.

A description of abelian extensions of real quadratic number fields in terms of coordinates of points of finite order on abelian varieties associated with certain modular curves was obtained by [Shimura  1972]   \cite{Shi1}.  
 Around that time,  [Stark 1976]   \cite{Sta1} formulated a number of conjectures on abelian extension of arbitrary number fields,  which in the real quadratic case amount to specifying generators of these extensions using special values of Artin $L$-functions.  In the case of an arbitrary number field, a description of the abelian extensions is given by class field theory, but an explicit description of the generators of these abelian extensions, in the sense sought by Kronecker and Hilbert, is still open.

The classical theory of complex multiplication on elliptic curves shows that the maximal abelian extension of 
${\Bbb Q}(\tau)$,  where $\tau$  is an imaginary quadratic irrationality,  can be obtained by adjoining the special 
values $\wp(\tau, z)$  and  $j(\tau)$  of modular functions $j$  and elliptic functions $\wp$,  and roots of unity, 
where $\tau$  generates the imaginary quadratic field and $z$  represents a torsion point on the corresponding elliptic curve. In particular, when $\tau$  is imaginary quadratic,  so that the lattice ${\Bbb Z}+{\Bbb Z}\tau$ has complex multiplication, the field  ${\Bbb Q}(\tau,  j(\tau))$  is the Hilbert class field of ${\Bbb Q}(\tau)$,  that is the maximal abelian unramified extension of ${\Bbb Q}(\tau)$.   In particular ${\Bbb Q}(\tau,  j(\tau))$  is a finite Galois extension of  ${\Bbb Q}(\tau)$  and $[{\Bbb Q}(\tau, j(\tau)) : {\Bbb Q}(\tau)]$ equals the class number of ${\Bbb Q}(\tau)$. 
Moreover,  the ideal class group of ${\Bbb Q}(\tau, j(\tau))$  is isomorphic to the Galois group of 
${\Bbb Q}(\tau, j(\tau))$  over ${\Bbb Q}(\tau)$,  whose action on $j(\tau)$  can be explicitly described,
see e.g. [Silverman 1994]   \cite[Ch. II]{S}.

In his theory of real multiplication  [Manin 2004]   \cite{Man1} proposed finding the analogue of this very specific type of result for real quadratic fields,  using so-called  ``pseudo-lattices'' 
${\Bbb Z}+{\Bbb Z}\theta$, where $\theta$  is real,   with non-trivial real multiplications.
In the classical theory,  a lattice is a ${\Bbb Z}$-module ${\Bbb Z}+{\Bbb Z}\tau$  of rank 2 with $\tau$ not real, 
isomorphic lattices differing by a non-zero scaling in ${\Bbb C}$.   The lattice ${\Bbb Z}+{\Bbb Z}\tau$ has non-trivial endomorphisms if and only if $\tau$  is imaginary quadratic,  and the endomorphism ring is then an 
order in the  field ${\Bbb Q}(\tau)$.  In Manin's theory,  a pseudo-lattice is a 
${\Bbb Z}$-module ${\Bbb Z}+{\Bbb Z}\theta$  of 
rank 2 with $\theta$  real.  Pseudo-lattices are again considered up to scaling,  but have non-trivial endomorphism 
rings if and only if $\theta$ is real quadratic.  The endomorphism ring $R={\Bbb Z}+\mathfrak{f} O$ is then an order 
in the field $k={\Bbb Q}(\theta)$,  where  $O$  is the ring of integers of $k$, and $\mathfrak{f}$  is the conductor of $R$. Manin calls these pseudo-lattices with real multiplication;  we shall call it here  {\it Manin pseudo-lattices}.

Let $\Lambda_M$ be Manin's  pseudo-lattice and  $R=End~(\Lambda_M)$ its endomorphism ring;  denote 
by $k=R\otimes {\Bbb Q}$ the real quadratic field associated to the ring $R$;   notice  that  $R={\Bbb Z}+\mathfrak{f}O_k$ is an order in $k$, 
where $\mathfrak{f}\ge1$ is the conductor of the order.  
The class field theory says that for each abelian extension $K$ of the field $k$ it holds
$Gal~(K|k)\cong Cl~(k)$, where $Cl~(k)$ is the (abelian) group of ideal classes of the ring $R$;
such an extension is called a {\it ring class field} of $k$  modulo  (conductor)  $\mathfrak{f}$. 
Notice that the ring class field modulo $\mathfrak{f}=1$  coincides  with the Hilbert class field,  i.e. the maximal 
unramified extension  of $k$ [Silverman 1994]   \cite{S}.  We shall focus on the following
 real multiplication problem: 
{\it  To  construct explicit generators of the ring class field of $k$ modulo  $\mathfrak{f}\ge 1$.}

\medskip
Let $N>0$ be an integer and $\Gamma_0(N):=\{(a,b,c,d)\in SL_2({\Bbb Z}) ~|~ c\equiv 0 ~mod~ N\}$ 
a finite index  subgroup of the modular group.
By $S_2(\Gamma_0(N))$ one understands a collection of the cusp forms
(of weight two) on the extended upper half-plane 
 ${\Bbb H}^*={\Bbb H}\cup {\Bbb Q}\cup\{\infty\}$. 
Recall that the Riemann surface $X_0(N):= {\Bbb H}^*/\Gamma_0(N)$ is called 
a {\it modular curve}.    We identify  $S_2(\Gamma_0(N))$ with the linear space 
$\Omega_{hol}(X_0(N))$ of the holomorphic differentials  on $X_0(N)$  via the formula 
$f(z)\mapsto \omega=f(z)dz$.  By ${\Bbb T}_{\Bbb Z}:= {\Bbb Z}[T_1, T_2,\dots]$
we denote a commutative algebra of  the Hecke operators $T_n$ acting on the space  
$S_2(\Gamma_0(N))$ by the formula $T_n f=\sum_{m\in {\Bbb Z}}\gamma(m)q^m$, where
$\gamma(m)= \sum_{a|GCD~(m,n)}a c_{mn/a^2}$ and  $f(z)=\sum_{m\in {\Bbb Z}}c(m)q^m$ is the Fourier
series of the cusp form $f$ at $q=e^{2\pi iz}$.
The common eigenvector $f\in S_2(\Gamma_0(N))$ of all $T_n\in {\Bbb T}_{\Bbb Z}$
is referred to as a {\it Hecke eigenform}.   By  $K_f$ we understand  an 
algebraic number field generated by the Fourier coefficients of the Hecke eigenform $f$;
it is  known that $K_f$ is  totally real and $deg~(K_f | {\Bbb Q})\le g$, where  $g$ is  the genus of surface  
$X_0(N)$ [Diamond \& Shurman 2005]   \cite{DS}, pp. 234-235.  When $deg~(K_f | {\Bbb Q})=g$ the Hecke eigenform $f$ will be called  {\it maximal};
in what follows we work exclusively with this class of the Hecke eigenforms.  
For a Hecke eigenform $f$ one  considers  a ${\Bbb Z}$-module
$\int_{H_1(X_0(N), Sing~\omega; {\Bbb Z})} \Re~(\omega) = {\Bbb Z}\lambda_1+\dots+{\Bbb Z}\lambda_g$,
where $\omega=f(z)dz$ is a holomorphic differential on $X_0(N)$; 
such a module belongs to the field $K_f$  (lemma \ref{lm1}).
The vector $v_A=(\lambda_1,\dots,\lambda_g)$  coincides with the Perron-Frobenius
eigenvector of a positive integer matrix $A$ corresponding to   the 
eigenvalue $\lambda_A\in K_f$;  the algebraic number $\lambda_A$ we 
shall call a {\it Hecke unit} of the field $K_f$.

On the other hand,  every measured foliation $\mathcal{F}$
on a surface $X$ 
([Hubbard \& Masur 1979]   \cite{HuMa1}, [Thurston 1988]   \cite{Thu1}) 
defines a  pseudo-lattice $\Lambda=Jac~(h(\mathcal{F}))$,
where $h(\mathcal{F})$ is an induced measured foliation on the torus and $Jac~(h(\mathcal{F}))$
its jacobian (lemma \ref{lm2});  in particular  if  $\mathcal{F}=\mathcal{F}_N$
is a measured foliation of the surface $X_0(N)$ by the vertical trajectories $\Re~(fdz)=0$ of a 
Hecke eigenform $f$ then $Jac~(h(\mathcal{F}_N))$ is a Manin pseudo-lattice $\Lambda_M$ (lemma \ref{lm5}).  
By a {\it j-invariant} of $\Lambda_M$ we understand 
the algebraic number $j(\Lambda_M):=\lambda_A$, where $\lambda_A$ is the Hecke
unit of the  field $K_f$;  it is indeed an invariant independent of  basis
in the ${\Bbb Z}$-modules  ${\Bbb Z}\lambda_1+\dots+{\Bbb Z}\lambda_g$
and  $Jac~(h(\mathcal{F}_N))$, and as such generalizes the classical j-invariant
defined on  the upper-half plane ${\Bbb H}$ to the (quadratic irrational) 
points at the boundary of ${\Bbb H}$.  We shall consider the ring 
$R={\Bbb Z}+\mathfrak{f}O_k:=End~(\Lambda_M)$ and  the real quadratic  
field $k=R\otimes {\Bbb Q}$.   Our main result can be formulated  as follows.
\begin{theorem}\label{thm1}
The extension $K=k(j(\Lambda_M))$ is a ring class field of $k$ modulo $\mathfrak{f}\ge 1$.
 \end{theorem}
The structure of the article is as follows. The notation and preliminary facts
are introduced in Section 2.  Theorem \ref{thm1} is proved in Section 3.
In Section 4 we compare theorem \ref{thm1}
with the results of  [Shimura  1972]   \cite{Shi1}.

\section{Preliminaries}
\subsection{Measured foliations and their jacobians}
By a $p$-dimensional  $C^r$ foliation of an $m$-dimensional 
manifold $M$ one understands a decomposition of $M$ into a union of
disjoint connected subsets $\{ \mathcal{L}_{\alpha}\}_{\alpha\in A}$  called
{\it leaves} of the foliation. The leaves must satisfy the
following property: Every point in $M$ has a neighborhood $U$
and a system of local class $C^r$ coordinates 
$x=(x^1,\dots, x^m): U\to {\Bbb R}^m$ such that for each leaf 
$\mathcal{L}_{\alpha}$, the components of $U\cap \mathcal{L}_{\alpha}$
are described by the equations $x^{p+1}=Const, \dots, x^m=Const$.
Such a foliation is denoted by $\mathcal{F}=\{ \mathcal{L}_{\alpha}\}_{\alpha\in A}$.
The number $q=m-p$ is called a {\it codimension} of the foliation
$\mathcal{F}$ [Lawson 1974]   \cite{Law1} p.370.    
The codimension $q$ class $C^r$ foliations $\mathcal{F}_0$ and $\mathcal{F}_1$ 
are said to be $C^s$-conjugate ($0\le s\le r$) if there exists a diffeomorphism
of $M$  of class $C^s$, which maps the leaves of $\mathcal{F}_0$
onto the leaves of $\mathcal{F}_1$.  If $s=0$  $\mathcal{F}_0$ and $\mathcal{F}_1$
are {\it topologically conjugate}, {\it ibid.}, p.388.

The foliation $\mathcal{F}$ is called {\it singular} if the codimension
$q$ of the foliation depends on the leaf.  We further  assume
that $q$ is constant for all  but a {\it finite}
number of leaves. Such a set of the exceptional leaves will be denoted by
$Sing ~\mathcal{F}:= \{ \mathcal{L}_{\alpha}\}_{\alpha\in E}$, where $|E|<\infty$;
note that if  $Sing ~\mathcal{F}=\emptyset$  one gets the usual definition of a 
(non-singular) foliation.  A quick  example of the singular foliations is given
by trajectories of a non-trivial differential form  on the manifold $M$
vanishing  in a finite number of points of $M$; 
the set of zeroes of such a form corresponds to the  exceptional leaves of 
foliation.

Roughly speaking,  measured foliation ([Hubbard \& Masur 1979]   \cite{HuMa1}, \cite[Section 0.3.2]{N})  is a singular codimension one  $C^r$-foliation,
induced  by  trajectories  of  a closed differential 
 $\phi$ on a two-dimensional manifold (surface) $X$;  formally speaking,  
 a {\it measured foliation}   $\mathcal{F}$ on a surface $X$
is a  partition of $X$ into the singular points $x_1,\dots,x_n$ of
order $k_1,\dots, k_n$ and regular leaves (1-dimensional submanifolds). 
On each  open cover $U_i$ of $X-\{x_1,\dots,x_n\}$ there exists a non-vanishing
real-valued closed 1-form $\phi_i$  such that: 
(i)  $\phi_i=\pm \phi_j$ on $U_i\cap U_j$;
(ii) at each $x_i$ there exists a local chart $(u,v):V\to {\Bbb R}^2$
such that for $z=u+iv$, it holds $\phi_i=Im~(z^{k_i\over 2}dz)$ on
$V\cap U_i$ for some branch of $z^{k_i\over 2}$.
The pair $(U_i,\phi_i)$ is called an atlas for the measured foliation $\mathcal{F}$.
Finally, a measure $\mu$ is assigned to each segment $(t_0,t)\in U_i$, which is  transverse to
the leaves of $\mathcal{F}$, via the integral $\mu(t_0,t)=\int_{t_0}^t\phi_i$. The 
measure is invariant along the leaves of $\mathcal{F}$, hence the name. 
Note that  measured foliation can have  singular points  of the half-integer index;  
those  cannot be given by the trajectories of a closed form. Yet
when all $k_i$ are even integers,   the leaves of foliation can be continuously 
oriented  and   $\mathcal{F}$ is called {\it oriented} in this case;   such foliations are given 
by trajectories of  a closed differential form on the surface $X$. 
In what follows  we  work with the  oriented measured foliations.

Let $\mathcal{F}$ be  measured foliation on a compact surface
$X$ and  $\{\gamma_1,\dots,\gamma_n\}$  a
basis in the relative homology group $H_1(X, Sing~\mathcal{F}; {\Bbb Z})$;
it is  known  that $n=2g+|Sing~(\mathcal{F})|-1$, where $g$ is the genus of $X$
[Hubbard \& Masur 1979]   \cite{HuMa1}.  By $\lambda_i\in {\Bbb R}$ we denote  
the periods $\int_{\gamma_i}\phi$ of $\phi$ in the above basis.
By a {\it jacobian} of the measured foliation $\mathcal{F}$ one 
understand a ${\Bbb Z}$-module ${\Bbb Z}\lambda_1+\dots+{\Bbb Z}\lambda_n$
regarded as a subset of the real line ${\Bbb R}$;  we shall denote the jacobian by  $Jac~(\mathcal{F})$.
The jacobian  is independent of the choice of 
basis in the homology group $H_1(X, Sing~\mathcal{F}; {\Bbb Z})$ and depends solely on the 
foliation $\mathcal{F}$ \cite{Nik1}; moreover the measured foliations $\mathcal{F}$ and $\mathcal{F}'$  
are topologically  conjugate if and only if $Jac~(\mathcal{F}')=\mu ~Jac~(\mathcal{F})$,
where $\mu>0$ is a real number {\it ibid}.

Let $\varphi: X\to X$ be an orientation-preserving automorphism of a compact surface $X$;
it is known that such an  automorphism is  either
(i) of finite order or (ii) an infinite order automorphism preserving certain  simple closed curves
or else (iii) a  pseudo-Anosov,  i.e. an  infinite order automorphism which  does not preserve any  
simple closed curve on $X$ [Thurston 1988]   \cite{Thu1}.  In case (iii)  there exist a stable $\mathcal{F}_s$
and unstable $\mathcal{F}_u$ mutually orthogonal measured foliations on $X$
such that $\varphi(\mathcal{F}_s)={1\over\lambda_{\varphi}}\mathcal{F}_s$ 
and $\varphi(\mathcal{F}_u)=\lambda_{\varphi}\mathcal{F}_u$,  where $\lambda_{\varphi}>1$
is called the {\it dilatation} of $\varphi$.  The invariant measured foliation $\mathcal{F}=\mathcal{F}_u$
is called a {\it pseudo-Anosov foliation};  it is known that its jacobian belongs to 
the number field ${\Bbb Q}(\lambda_{\varphi})$ [Thurston 1988]   \cite{Thu1}, p.427-428.

\subsection{Foliations on modular curves}
Let $N>1$ be a natural number and  consider a finite index subgroup 
of the modular group given by the formula
$\Gamma_0(N) = \{(a, b,  c, d)\in SL_2({\Bbb Z})~|~
c\equiv 0~mod ~N\}$.
Let ${\Bbb H}=\{z=x+iy\in {\Bbb C} ~|~ y>0\}$ be the upper half-plane  and 
let $\Gamma_0(N)$  act on ${\Bbb H}$  by the linear fractional
transformations;  consider an orbifold  ${\Bbb H}/\Gamma_0(N)$.
To compactify the orbifold 
at the cusps, one adds a boundary to ${\Bbb H}$,  so that 
${\Bbb H}^*={\Bbb H}\cup {\Bbb Q}\cup\{\infty\}$ and the compact Riemann surface 
$X_0(N)={\Bbb H}^*/\Gamma_0(N)$ is called a {\it modular curve}.   
The meromorphic functions $f(z)$ on ${\Bbb H}$ that
vanish at the cusps and such that
$f\left({az+b\over cz+d}\right)=(cz+d)^2f(z),
~\forall \left(\small\begin{matrix} a & b\cr c & d\end{matrix}\right)\in\Gamma_0(N)$, 
are  called  {\it cusp forms} of weight two;  the (complex linear) space of such forms
will be denoted by $S_2(\Gamma_0(N))$.  The formula $f(z)\mapsto \omega=f(z)dz$ 
defines an isomorphism  $S_2(\Gamma_0(N))\cong \Omega_{hol}(X_0(N))$, where 
$\Omega_{hol}(X_0(N))$ is the space of holomorphic differentials
on the Riemann surface $X_0(N)$.  Note that 
$\dim_{\Bbb C}(S_2(\Gamma_0(N))=\dim_{\Bbb C}(\Omega_{hol}(X_0(N))=g$,
where $g=g(N)$ is the genus of the surface $X_0(N)$. 
Recall that there exists a natural involution $i$ on the space $S_2(\Gamma_0(N))$
defined by the formula $f(z)\mapsto f^*(z)$, where $f(z)=\sum c_nq^n$ and 
$f^*(z)=\sum \bar c_nq^n$. A subspace, $S_2^{\Bbb R}(\Gamma_0(N))$, fixed by the
involution, consists of the cusp forms, whose Fourier coefficients are
the real numbers. Clearly,  $dim_{\Bbb R}(S_2^{\Bbb R}(\Gamma_0(N)))=g$.
The map $i$ induces an involution $i_{\Phi}$ on the space $\Phi_{X_0(N)}$ of all measured foliations on $X_0(N)$;
in a proper coordinate system the involution $i_{\Phi}$  acts by the
formula $(\lambda_1,\dots,\lambda_g,\lambda_1',\dots,\lambda_g')\mapsto
(\lambda_1',\dots,\lambda_g',\lambda_1,\dots,\lambda_g)$. 
A  subspace of $\Phi_{X_0(N))}$ fixed by involution $i_{\Phi}$ we shall denote by  $\Phi_{X_0(N)}^{\Bbb R}$; 
it consists of  measured foliations of the form  $(\lambda_1,\dots,\lambda_g,\lambda_1,\dots,\lambda_g)$.
Thus  $Jac~(\mathcal{F})={\Bbb Z}\lambda_1+\dots+{\Bbb Z}\lambda_g$ for 
$\forall \mathcal{F}\in \Phi_{X_0(N)}^{\Bbb R}$.

Let $f\in S_2(\Gamma_0(N))$ be a (normalized) Hecke eigenform, such that
\linebreak
$f(z)=\sum_{n=1}^{\infty}c_n(f)q^n$ its Fourier series. We shall denote
by $K_f={\Bbb Q}(\{c_n(f)\})$ the algebraic number field generated
by the Fourier coefficients of $f$. Let $g$ be the genus of the modular curve $X_0(N)$.
It is well known that $1\le deg~(K_f~|~{\Bbb Q})\le g$
and $K_f$ is a totally real field, see e.g. [Diamond \& Shurman 2005]   \cite[pp. 234-235]{DS}.
Let $\mathcal{F}_N$ be a measured foliation given by the lines $Re~(fdz)=0$,  where $f$ is the maximal
Hecke eigenform.  We shall use   the following fact.  
\begin{lemma}\label{lm1}
{\bf (\cite{Nik1})}
$\mathcal{F}_N$ is pseudo-Anosov and 
$Jac~(\mathcal{F}_N)$ is a ${\Bbb Z}$-module in the field $K_f$. 
\end{lemma}

\section{Proof of theorem 1}
For the sake of clarity, we  outline  main ideas of the proof. 
Given a measured foliation $\mathcal{F}$ on a surface $X$, 
one can associate to $\mathcal{F}$ the canonical measured 
foliation, $F$, on torus $T^2$ and {\it vice versa} provided
the set $Sing~\mathcal{F}$ is specified.  This is the standard 
fact following from the Riemann-Hurwitz construction of a
ramified covering of one Riemann surface (torus $T^2$) 
by another (surface $X$).    Such a map between foliations 
is well defined and unique up to the topological conjugacy.

Since $Jac~(\mathcal{F})$ is  invariant of the topological conjugacy
classes,  $Jac~(\mathcal{F})$ must be  related to $Jac~(F)$;
it is indeed so and such a relation is the inclusion $Jac~(F)\subseteq Jac~(\mathcal{F})$.

In case $\mathcal{F}$ is a pseudo-Anosov measured foliation,   the
inclusion $Jac~(F)\subseteq Jac~(\mathcal{F})$ can be exploited
to construct (finite)  extensions of the real quadratic number fields.   
Namely,  whenever $\mathcal{F}$ is such a foliation,  it is known that
$Jac~(\mathcal{F})$ belongs (as a ${\Bbb Z}$-module) to a real algebraic
field $K$ generated by dilatation $\lambda_{\varphi}$ of
the pseudo-Anosov automorphism $\varphi$ of surface $X$. 
The automorphism descends canonically to such of torus $T^2$
and,  therefore,  $Jac~(F)$ belongs to a real quadratic field $k$;
in such a way one gets an extension $k\subseteq K$ of the 
field $k$.

The above construction is general   and  (in most cases) the arithmetic
   of the  field $k$ and such of the field  $K$  are  independent of each other.
However,  it is not so for the modular curve $X=X_0(N)$ 
and the foliation $\mathcal{F}_N$  given by the vertical trajectories of 
the corresponding Hecke eigenform.  It will develop,   that in this  special case 
the ideal class group $Cl~(k)$ of the field $k$ 
controls    the Galois group  $Gal~(K|k)$  of extension $k\subseteq K$,
so that   $Gal~(K|k)\cong Cl~(k)$.

Why  $Gal~(K|k)\cong Cl~(k)$ anyway?   
Recall that there are $\{f_1,\dots, f_g\}$ linearly independent Hecke 
eigenforms,  where $g$ is the genus of $X_0(N)$.   
For each $1\le i\le g$ the Fourier coefficients of eigenform $f_i$ 
generate an algebraic number field $K=K_f$ such that
$deg~(K_f|{\Bbb Q})=g$  and $Jac~(\mathcal{F}_i)\subset K$,
where $\mathcal{F}_i$  is a measured foliation by  vertical 
trajectories of the form $f_idz$. 
On the other hand,  it was shown that on the torus $T^2$ 
one gets  $\{F_1,\dots, F_g\}$ measured foliations corresponding
to $\mathcal{F}_i$.   Because $f_i$ are  conjugate Hecke eigenforms, 
the  jacobians $Jac~(F_i)$  will have the same  endomorphism ring,
yet  being  itself  pairwise distinct;  it means that $|Cl~(k)|=g$. 
With  a little extra work (see below), one constructs the 
required isomorphism $Gal~(K|k)\cong Cl~(k)$.

In view of the above construction, $Jac~(F_i)$ are Manin 
pseudo-lattices;  therefore one gets a solution to the real
multiplication problem of Yu.~I.~Manin. 
(The reader can verify that all maps above are equivariant 
and canonical with respect to the natural morphisms between
the objects in use.)   These remarks complete the sketch of proof of our 
main result.

\bigskip
We shall pass to a detailed argument by 
splittting  the proof in a series of lemmas starting with the
following Riemann-Hurwitz construction for  measured foliations. 
\begin{lemma}\label{lm2}
Every  measured foliation $\mathcal{F}$  on a surface $X$  (with fixed set $Sing~\mathcal{F}$)  covers
a measured foliation $h(\mathcal{F})$ on torus $T^2$ and vice versa;  
the map $h$ is a bijection between the classes of topologically conjugate 
measured foliations on $X$ and such on $T^2$.    
\end{lemma}
\begin{proof}
Let $\mathcal{F}$ be measured foliation on a surface $X$; let
$Jac~(\mathcal{F})={\Bbb Z}\lambda_1+\dots+{\Bbb Z}\lambda_n$
be its jacobian. We shall consider an $n$-dimensional torus
$T^n={\Bbb R}^n/{\Bbb Z}^n$ endowed with a codimension one
foliation $\mathcal{F}'$ defined by the formula
$\int_{H_1(T^n; {\Bbb Z})}\phi_{\mathcal{F}'}=Jac~(\mathcal{F})$,
where $\phi_{\mathcal{F}'}$ is a closed differential tangent to  
 foliation $\mathcal{F}'$;  such a foliation always exists [Plante  1975]   \cite{Pla1}.
Notice that $Jac~\mathcal{F}$ coincides with the so-called {\it group
of periods}  $P(\mathcal{F}'):=\int_{H_1(T^n; {\Bbb Z})}\phi_{\mathcal{F}'}$
of foliation $\mathcal{F}'$; such a group  is  defined in [Plante  1975]   \cite{Pla1},  p.346.

Let $T^2$ be a two-dimensional torus; let $E=(T^n, T^2, p)$ be
a fiber bundle $T^n$ over $T^2$ with the fiber $p^{-1}(x)\cong T^{n-2}$
over the point $x\in T^2$.  Recall that every foliation $F$ 
on $T^2$  induces a foliation $p^{-1}(F)$ on $T^n$
defined by the fiber map $p: T^n\to T^2$ [Lawson 1974]   \cite{Law1}, p. 373.
Denote by $F$ a foliation on the torus $T^2$ such that
$p^{-1}(F)=\mathcal{F}'$. Since $\mathcal{F}'$ is given by a closed form, so will be the
foliation $F$; in other words, $F$ is a measured foliation. 
(Notice  that the rank of $Jac~(F)$ is less or equal to such of
$Jac~(\mathcal{F}')$  in  general;   however,  $rank~(Jac~(F))=2$
and  $rank~(Jac~(\mathcal{F}'))=n$  for generic foliations  $F$ and $\mathcal{F}'$.)   
Thus we have  a correctly defined map
\begin{equation}
h: \mathcal{F}\longrightarrow \mathcal{F}'\longrightarrow F. 
\end{equation}
If $E'=(T^n, T^2, p')$ is another
fiber bundle with the base $T^2$ then the corresponding foliation
$F'$ will be topologically conjugate to the foliation $F$;  we
leave this claim as an exercise to the reader.

Conversely,  let $F$ be measured foliation on the
two-torus $T^2$; if $E=(T^n,T^2,p)$ is a fiber bundle,
then $\mathcal{F}'=p^{-1}(F)$ is an induced codimension one 
foliation on $T^n$ [Lawson 1974]   \cite{Law1}, p. 373.  The latter is 
given by a closed form $\phi_{\mathcal{F}'}$ and has the
group of periods $P(\mathcal{F}')= \int_{H_1(T^n; {\Bbb Z})}\phi_{\mathcal{F}'}$ 
[Plante  1975]   \cite{Pla1}, p. 346.  Fixing a set $Sing~\mathcal{F}$ 
one constructs measured foliation $\mathcal{F}$ on a surface $X$,
such that $Jac~(\mathcal{F})=P(\mathcal{F}')$. Clearly, every 
measured  foliation $\mathcal{F}$ on $X$ can be constructed in this
way. Thus the map $h: \mathcal{F}\to F$ is a bijection between the
classes of topological conjugacy  of 
measured foliations on $X$ (with fixed set $Sing~\mathcal{F}$) and
such on the torus $T^2$. Lemma \ref{lm2} follows.
\end{proof}

\begin{lemma}\label{lm3}
If $\mathcal{F}$ is  measured foliation  on a surface $X$ then
$Jac~(h(\mathcal{F}))\subseteq Jac~(\mathcal{F})$.  
\end{lemma}
\begin{proof}
Recall that the integral $\int_{\gamma}\phi$ defines 
a pairing $H_1(M; {\Bbb R})\times H^1(M; {\Bbb R})\rightarrow {\Bbb R}$
between the first homology and cohomology groups of a manifold $M$.
The pairing is a scalar product on the vector space $H_1(M; {\Bbb R})
\cong H^1(M; {\Bbb R})$; we shall fix an orthonormal basis $(e_i,e_j)=\delta_{ij}$
in $H_1(M; {\Bbb R})$. 
Let $X$ be a surface with measured foliation $\mathcal{F}$. It was shown
that $Jac~(\mathcal{F})=Jac~(\mathcal{F}')$, where $\mathcal{F}'$ is a codimension
one foliation on the torus $T^n$ and $F=h(\mathcal{F})$ is a foliation on $T^2$
induced by the fiber map $p: T^n\to T^2$. Since $p$ is a continuous map
it defines a homomorphism
$p_*: H_1(T^n; {\Bbb Z})\cong {\Bbb Z}^n \to  H_1(T^2; {\Bbb Z})\cong {\Bbb Z}^2$.
In a proper basis  the map $p_*$ coincides with the projection map  on the first 
coordinates, i.e. a map  $p_*: {\Bbb Z}^n\to {\Bbb Z}^2$ which acts  by
the formula $(z_1,z_2,\dots, z_n)\mapsto (z_1,z_2, 0,\dots,0)$.   
One can compare the following system of equations
\begin{equation}
\left\{
\begin{array}{ccc}
Jac~(\mathcal{F}) &= \int_{H_1(T^n; {\Bbb Z})}\phi &=\lambda_1 {\Bbb Z}+\dots+\lambda_n {\Bbb Z}\\
Jac~(h(\mathcal{F})) &= \int_{H_1(T^2; {\Bbb Z})}p^*(\phi) &= \lambda_1 {\Bbb Z}+\lambda_2 {\Bbb Z},
\end{array}
\right.
\end{equation}
where $\lambda_i$ are coordinates of the form $\phi$ in $H^1(T^n; {\Bbb R})$
and $p^*(\phi)$ its projection in $H^1(T^2; {\Bbb R})$.  Thus 
$Jac~(h(\mathcal{F}))\subseteq Jac~(\mathcal{F})$ and the equality holds 
if and only if $n=2$.  Lemma \ref{lm3} is proved. 
\end{proof}

\begin{lemma}\label{lm4}
If $\mathcal{F}$ is a pseudo-Anosov  measured foliation on a surface $X$  then 
$K=End~(Jac~(\mathcal{F}))\otimes ~{\Bbb Q}$   is a number field;
the field $K=k(\lambda_{\varphi})$, where $k=End~(Jac~(h(\mathcal{F})))\otimes ~{\Bbb Q}$ 
is either a real quadratic or rational field.
\end{lemma}
\begin{proof}
Let $\mathcal{F}$ be measured foliation which is invariant of an
infinite-order (pseudo-Anosov) automorphism $\varphi$ of a
surface $X$; let $\lambda_{\varphi}>1$ be the dilatation of
automorphism $\varphi$. It is known that in this case the generators 
$(\lambda_1,\dots,\lambda_n)$ of $Jac~(\mathcal{F})$ are coordinates of
the Perron-Frobenius eigenvector corresponding to the eigenvalue $\lambda_{\varphi}$
of a positive integer matrix $A_{\varphi}\in GL_n({\Bbb Z})$. 
Therefore $\lambda_i$ can be scaled to belong to the algebraic number
field $K={\Bbb Q}(\lambda_{\varphi})$; clearly the endomorphism ring
$End~(Jac~(\mathcal{F}))$ is an order in the field $K$ and $K$ itself 
can be written as   $K=End~(Jac~(\mathcal{F}))\otimes ~{\Bbb Q}$.

The matrix $A_{\varphi}$ defines an automorphism of the torus $T^n$;  
the automorphism commutes with the projection map 
$p_*:H_1(T^n; {\Bbb Z})\to H_1(T^2; {\Bbb Z})$ defined earlier.
We shall write $A\in GL_2({\Bbb Z})$ to denote an automorphism
of $T^2$ induced by $A_{\varphi}$. Since $A_{\varphi}$ has an
infinite order so will be the automorphism $A$;  therefore
the eigenvalues $\lambda$ and $\bar\lambda$ of $A$ are either 
quadratic irrational numbers or $\lambda=\bar\lambda=\pm 1$.
When $\lambda>1$ is a quadratic irrationality the $Jac~(h(\mathcal{F}))$
is generated by coordinates of the eigenvector corresponding to $\lambda$;
the latter can be scaled to belong to the field $k={\Bbb Q}(\lambda)$.
Thus   $End~(Jac~(h(\mathcal{F})))\otimes ~{\Bbb Q}$ is a real quadratic
number field which coincides with the field $k$.
When  $\lambda=\bar\lambda=\pm 1$ then $Jac~(h(\mathcal{F}))$ belongs to 
the field ${\Bbb Q}$ so that  $End~(Jac~(h(\mathcal{F})))\otimes ~{\Bbb Q}\cong {\Bbb Q}$.
In view of lemma \ref{lm3}  one arrives at the following inclusion 
 of the number fields
\begin{equation}
End~(Jac~(h(\mathcal{F})))\otimes ~{\Bbb Q}\subseteq End~(Jac~(\mathcal{F}))\otimes ~{\Bbb Q}.
\end{equation}
In other words, the field $K=k(\lambda_{\varphi})$,  where $k$ 
is either a real quadratic or rational number field.
Lemma \ref{lm4} follows.
\end{proof}

\begin{remark}
Notice that in general the extension $K|k$ is not a Galois 
extension; if it is such an extension, the group $Gal~(K|k)$
is not necessarily defined by the ideal class group $Cl~(k)$.
Yet $Gal~(K|k)\cong Cl~(k)$ whenever $\mathcal{F}=\mathcal{F}_N$
is measured foliation by the vertical trajectories of a Hecke
eigenform $f\in S_2(\Gamma_0(N))$ (see lemma \ref{lm6}).   
\end{remark}
Recall that $f\in S_2(\Gamma_0(N))$ is called a maximal Hecke
eigenform if $deg~(K_f)=g$,  where $g$ is the genus of surface 
$X_0(N)$; consider a measured foliation $\mathcal{F}_N=\Re~(fdz)$
generated by $f$. 
\begin{lemma}\label{lm5}
$Jac~(h(\mathcal{F}_N))$ is  a Manin  pseudo-lattice.
\end{lemma}
\begin{proof}
The foliation $\mathcal{F}_N$ is the invariant measured foliation of
a pseudo-Anosov automorphism of $X_0(N)$ (lemma \ref{lm1});
thus to prove our claim it remains to show that 
$End~(Jac~(h(\mathcal{F}_N)))\otimes ~{\Bbb Q}\not\cong {\Bbb Q}$ (lemma \ref{lm4}). 
Consider a Hecke operator  $T_n\in {\Bbb T}_{\Bbb Z}$; it acts on 
$Jac ~(\mathcal{F}_N)={\Bbb Z}\lambda_1+\dots+{\Bbb Z}\lambda_g$ by 
a symmetric matrix $(t_{ij})\in M_g({\Bbb Z})$, i.e. $\lambda_j'=\sum t_{ij}\lambda_i$
[Diamond \& Shurman 2005]   \cite{DS}.   Let us calculate an induced action of $T_n$ on   
$Jac ~(h(\mathcal{F}_N))={\Bbb Z}\lambda_1+{\Bbb Z}\lambda_2$.
In  notation of lemma \ref{lm5}  we have 
$p_*(\lambda_1,\dots,\lambda_g)=(\lambda_1, \lambda_2, 0,\dots, 0)$,
where $p_*$ is a projection map $H_1(X_0(N), Sing~(\mathcal{F}_N; {\Bbb Z})
\to H_1(T^2; {\Bbb Z})$; thus the action of $T_n$ on $Jac ~(h(\mathcal{F}_N))$
is given by a symmetric matrix $\tau_n=(t_{11},t_{12},t_{12},t_{22})\in M_2({\Bbb Z})$. 
Since  $\tau_n$ is an integer matrix it  defines  a non-trivial   endomorphism of $Jac ~(h(\mathcal{F}_N))$; 
such an endomorphism can  be given as multiplication of points of ${\Bbb R}$
by a real number $\alpha$, i.e.   $\alpha\lambda_1 = t_{11}\lambda_1+t_{12}\lambda_2, 
~\alpha\lambda_2 = t_{12}\lambda_1+t_{22}\lambda_2$.
Since  $\theta=\lambda_2/\lambda_1$, it satisfies  the equation
${\theta}={t_{12}+t_{22}\theta\over t_{11}+t_{12}\theta}$;
the latter  is equivalent to the quadratic equation 
$t_{12}\theta^2+(t_{11}-t_{22})\theta-t_{12}=0$
whose  determinant   $D=(t_{11}-t_{22})^2+4t_{12}^2\ge 0$.
Such a  quadratic equation has two real roots;  these roots cannot be rational 
since $\tau_n$ is a non-trivial endomorphism.  
Therefore  $\theta$ is a quadratic irrationality.
Lemma \ref{lm5}  follows.
\end{proof}

\begin{lemma}\label{lm6}
The number field $K=End~(Jac~(\mathcal{F}_N))\otimes ~{\Bbb Q}$ is 
a ring class field of the real quadratic field  $k=End~(Jac~(h(\mathcal{F}_N)))\otimes ~{\Bbb Q}$.
\end{lemma}
\begin{proof}
To outline the proof, 
let  $R=End~(Jac~(h(\mathcal{F}_N)))$ be an order in the real quadratic field
$k=R\otimes {\Bbb Q}$ and $Cl~(R)$ the ideal class group of the order $R$.
To prove  lemma \ref{lm6} one needs to verify  the isomorphism  $Cl~(R)\cong Gal~(K|k)$,
where $Gal~(K|k)$ is the Galois group of the extension $K|k$;   in turn, such an isomorphism 
can  be given by  the action of group $Cl~(R)$ on  generators of the
extension $K|k$.   To give a rough idea,  look  at a basis of the Hecke eigenforms $\{f_1,\dots,f_g\}$
in the space $S_2(\Gamma_0(N))$;  such a basis is known to exist, see e.g. [Diamond \& Shurman 2005]   \cite[Theorem 5.8.2]{DS}.  
Each $f_i$ defines  a Manin pseudo-lattice  and  one gets  $\{\Lambda_M^{(1)},\dots, \Lambda_M^{(g)}\}$
pairwise non-isomorphic pseudo-lattices;  notice  that  $End~(\Lambda_M^{(i)})\cong R$ for all $1\le i\le g$,
i.e. $h_R=g$, where $h_R$ is the class number of the order $R$.
But $f_i$ define the Hecke units  $\lambda_i$ of the field $K$ which are algebraically conjugate numbers;
the formula  $j(\Lambda_M^{(i)})=\lambda_i$  gives the required  action of the group $Cl~(R)$ on
the generators of the field $K$.  We shall  pass to a detailed argument.

\medskip
(i) For a Hecke eigenform $f$ consider its Hecke unit $\lambda\in K_f$;
denote by $p\in {\Bbb Z}[x]$ a minimal polynomial of the algebraic number
$\lambda=\lambda_1$. Since $p(x)$ splits in the totally real field $K_f$ one can
write 
\begin{equation}\label{eq28}
p(x)=(x-\lambda_1) \dots (x-\lambda_g),
\end{equation}
where $\lambda_i\in K_f$ are the Hecke units of some Hecke eigenforms
$f_i$;  the eigenforms make a basis $\{f_1,\dots,f_g\}$ of the space $S_2(\Gamma_0(N))$.
By lemma \ref{lm5} there exist Manin's pseudo-lattices $\Lambda_M^{(i)}$
corresponding to $f_i$ and we shall define the j-invariant of the latter
as  $j(\Lambda_M^{(i)}):=\lambda_i$.

\medskip
(ii) Let $R=End~(\Lambda_M^{(1)})$ and let $h_R$ be the class number
of the order $R$.  Then $h_R=g$,  where $g$ is the genus of
the surface $X_0(N)$.  

\smallskip\hskip0.5cm
(a) Let us show that $g\le h_R$.  Indeed, a  basis $\{f_1,\dots,f_g\}$
of the Hecke eigenforms $f_i$ in the space $S_2(\Gamma_0(N))$ 
defines Manin's  pseudo-lattices $\Lambda_M^{(1)},\dots, \Lambda_M^{(g)}$,
such that $End~(\Lambda_M^{(i)})=R$.  Thus, the order $R$ has at least $g$ ideal
classes, i.e. $g\le h_R$.

\smallskip\hskip0.5cm
(b) Let us show that $g\ge h_R$.  Indeed, let $\Lambda_M^{(1)},\dots, \Lambda_M^{h_R}$
be a full list of Manin's   pseudo-lattices in the order $R$.  Since $R\cong {\Bbb T}_{\Bbb Z}/{\Bbb I}$,
where ${\Bbb I}$ is a fixed ideal of  the ring of Hecke operators ${\Bbb T}_{\Bbb Z}$, we conclude that
there exists at least $h_R$  Hecke eigenforms in the 
space $S_2(\Gamma_0(N))$.  Thus $g\ge h_R$.

\medskip
It follows from (a) and (b) that $g=h_R$.

\medskip
(iii) Let us  establish  an explicit formula for the isomorphism 
$Cl~(R)\to Gal~(K|k)$;  since $Gal~(K|k)$  is an automorphism group of the
field $K$, it will suffice  to define  the action of an element $a\in Cl~(R)$
on  generators of $K$. 
Let $\Lambda_M^{(i)}\subseteq R$ be a Manin pseudo-lattice  and let
$[\Lambda_M^{(i)}]$ be its  ideal class in $R$. Since
$[\Lambda_M^{(i)}]\in Cl~(R)$,  the element $a\ast [\Lambda_M^{(i)}]\in Cl~(R)$
for all $a\in Cl~(R)$.  We let $\Lambda_M^{(j)}$ be a Manin pseudo-lattice,
such that $[\Lambda_M^{(j)}]=a\ast [\Lambda_M^{(i)}]$.  For the sake of brevity,
we  write $\Lambda_M^{(j)}=a\ast [\Lambda_M^{(i)}]$. The action of an 
element $a\in Cl~(R)$ on the generators $j(\Lambda_M^{(i)})$ of the field
$K$ is given by the following formula:
\begin{equation}\label{eq31}
a\ast j(\Lambda_M^{(i)}):= j(a\ast [\Lambda_M^{(i)}]), \qquad\forall a\in Cl~(R).
\end{equation}
We  leave it to the reader to verify that the last formula gives an isomorphism
$Cl~(R)\to Gal~(K|k)$;  this argument completes  the proof of lemma \ref{lm6}.
\end{proof}

\medskip
Theorem \ref{thm1} follows from lemma \ref{lm6}.
$\square$

\section{Shimura's method revisited}
Theorem \ref{thm1} says that there exists an integer $N_\mathfrak{f}(d)$,  such
that the Fourier coefficients of a Hecke eigenform $f\in S_2(\Gamma_0(N_\mathfrak{f}(d)))$
generate  the ring class field of a  real quadratic field $k={\Bbb Q}(\sqrt{d})$  modulo the conductor 
$\mathfrak{f}\ge 1$. In 1972  Goro Shimura calculated the integer $N_\mathfrak{f}(d)$ for certain prime values 
of $d$  [Shimura  1972]   \cite{Shi1}.  (The general formula can be found in \cite{Nik2}.)
In this section we use Shimura's results to illustrate theorem \ref{thm1} by proving the 
following proposition. 
\begin{proposition}\label{lm7}
If $p$ is a prime number such that $p\equiv 1 \mod 4$, 
then $N_\mathfrak{f}(p)=p$ for some  $\mathfrak{f}\ge 1$. 
The  corresponding values of the conductor $\mathfrak{f}$ and the $j$-invariant $j(\Lambda_M)$
for $p\le 197$   are compiled  in Table 1. 
\end{proposition}
\begin{proof}
Recall that $\Gamma_1(N)$ is a subgroup of $\Gamma_0(N)$ defined by the 
congruence relations $a\equiv d\equiv 1 ~mod ~N$,  see Section 2.2 for the 
notation;  consider a character $\psi$ of $({\Bbb Z}/N {\Bbb Z})^{\times}$ such 
that $\psi(-1)=1$.  
Let $S_2(\Gamma_1(N), \psi)$ be a subspace of the space of cusp forms $S_2(\Gamma_1(N))$
consisting of function $f$ such that   $f\left({az+b\over cz+d}\right)=\psi(d) (cz+d)^2f(z),
~\forall \left(
\begin{smallmatrix} a & b\cr c & d\end{smallmatrix}
\right)\in\Gamma_0(N)$.  
For each positive integer $n$ one can define the Hecke operator $T(n)$ on $S_2(\Gamma_1(N))$
 [Diamond \& Shurman 2005]   \cite[p. 168]{DS};  a restriction of $T(n)$ to $S_2(\Gamma_1(N),  \psi)$ will be 
denoted by $T(n, \psi)$.  The corresponding Hecke eigenform we shall write as
$f(z, \psi)$ and the algebraic number field generated by the Fourier coefficients 
of $f(z, \psi)$ we denote by $K_{f, \psi}$.  Unless $\psi$ is a trivial character,
$K_{f, \psi}$ is a CM-field,  i.e. a totally imaginary quadratic extension of a totally
real algebraic number field [Shimura  1972]   \cite[Proposition 1.3]{Shi1}.     
The maximal real subfield of the field $K_{f,\psi}$ will be denoted by $F$. 

We further assume that $\psi$ is a non-trivial real quadratic character  $\psi^2=1$.  
It is not hard to see, that such a character generates an automorphism of 
order 2 of the field $F$ and,  therefore,  $Gal~(K|{\Bbb Q})$ contains a subgroup of order 2.
 Let $k$ be the corresponding real quadratic subfield of $F$.  Shimura noticed that
$F$ can be interpreted as a ray class field of $k$ connected to an integral ideal $\mathfrak{c}$
of $F$.   In particular, if $N$ is a prime number $p\equiv 1~mod~4$ and $\psi(a)=\left({a\over N}\right)$,
then $k={\Bbb Q}(\sqrt{N})$ [Shimura  1972]   \cite[p.147]{Shi1}.  In this case  Shimura gives numerical examples for
levels $29\le N\le 233$ with conductors  $N(\mathfrak{c})$,  where $N(\mathfrak{c})$ is the norm
of ideal $\mathfrak{c}$ [Shimura  1972]   \cite[pp.149-150]{Shi1}.

Shimura's method is linked to the measured foliations as follows.  Let $\mathcal{F}_{\psi}:=\Re ~(f(z, \psi)dz)$ be 
measured  foliation on the surface $X_1(N)$.   Consider the jacobian  $Jac~(\mathcal{F}_{\psi})={\Bbb Z}\lambda_1+
\dots+{\Bbb Z}\lambda_n$ of foliation $\mathcal{F}_{\psi}$;  it is clear that  $\lambda_i\in F$. 
 We claim that if $\psi$ is a quadratic character and companions
 [Shimura  1972]   \cite[p.134]{Shi1}
of $f(z, \psi)$ span   $S_2(\Gamma_1(N),  \psi)$,   then 
\begin{equation}\label{eq6}
Jac~(\mathcal{F}_{\psi})=Jac~(\mathcal{F}),
\end{equation}
where $\mathcal{F}=\Re~(g(z)dz)$ is measured foliation on surface $X_0(N)$ 
corresponding to   a maximal Hecke eigenform $g\in S_2(\Gamma_0(N))$. 
Indeed,  we have $Jac~(\mathcal{F})\subseteq Jac~(\mathcal{F}_{\psi})$,   since 
$S_2(\Gamma_0(N))\subset S_2(\Gamma_1(N))$ and the Hecke operator
$T(n,\psi)$ acts on both spaces [Diamond \& Shurman 2005]   \cite[p.168]{DS}.  
Notice that for the trivial character $\psi=1$ the following three fields 
coincide:  $K_{f, 1}=F=K_g$.  Thus,  for quadratic character $\psi$ field $K_{f, \psi}$
is  a quadratic  extension of $F=K_g$.   
But elements of the field $K_{f,\psi}$ are  periods of the holomorphic
differential $f(z, \psi)dz$ on surface $X_1(N)$;  the real parts of such periods generate
the maximal real subfield of $K_{f,\psi}$,  i.e. the field $F$.  Thus the ranks of 
$Jac~(\mathcal{F}_{\psi})$ and  $Jac~(\mathcal{F})$ must coincide, hence equation (\ref{eq6}).

In view of formula (\ref{eq6}),  one concludes  that  Shimura's real quadratic field $k$ (modulo integral
ideal $\mathfrak{c}$)  attached  to  quadratic character $\psi$  coincides with such obtained  
from Manin's   pseudo-lattice $\Lambda_M$  modulo the conductor $\mathfrak{f}=N(\mathfrak{c})$. 
Table 1 comprises all Shimura's  examples except for the non-companionate cusp forms $N=109, 157$ 
and $229$.   In case $N=233$ the value of conductor $N(\mathfrak{c})$ is unknown,  hence must be
excluded.    Proposition \ref{lm7} follows.
\end{proof}

\begin{table}[h]
\begin{tabular}{c|c|c|c}
\hline
&&&\\
$p$ & $k$ & $\mathfrak{f}$ & $j(\Lambda_M)$\\
&&&\\
\hline
$29$ & ${\Bbb Q}(\sqrt{29})$ & $5$ & $1$\\
\hline
$37$ & ${\Bbb Q}(\sqrt{37})$ & $1$ & $1$\\
\hline
$41$ & ${\Bbb Q}(\sqrt{41})$ & $1$ & $1$\\
\hline
$53$ & ${\Bbb Q}(\sqrt{53})$ & $7$ & $1+\sqrt{2}$\\
\hline
$61$ & ${\Bbb Q}(\sqrt{61})$ & $13$ & $2+\sqrt{3}$\\
\hline
$73$ & ${\Bbb Q}(\sqrt{79})$ & $89$ & ${1\over 2}(1+\sqrt{5})$\\
\hline
$89$ & ${\Bbb Q}(\sqrt{89})$ & $5$ &a fundamental unit of the splitting field\\
           &        &                                       &of polynomial  $x^3+17x^2+83x+125=0$\\
\hline
$97$ & ${\Bbb Q}(\sqrt{97})$ & $467$ & $x^3+27x^2+204x+467=0$\\
 \hline
$101$ & ${\Bbb Q}(\sqrt{101})$ & $5$ & $x^4+13x^3+51x^2+67x+20=0$\\
 \hline
$113$ & ${\Bbb Q}(\sqrt{113}) $  & $97$ & $x^4+19x^3+122x^2+297x+194=0$\\
 \hline
$137$ & ${\Bbb Q}(\sqrt{137}) $ & $109$  & $x^5+23x^4+188x^3+670x^2+989x+436=0$\\
 \hline
$149$ & ${\Bbb Q}(\sqrt{149})$ & $61$ & $x^6+20x^5+148x^4+499x^3+766x^2+$\\
            &           &                                          & $+465x+61=0$\\
 \hline
$173$ & ${\Bbb Q}(\sqrt{173})$  & $13$ & $x^7+20x^6+151x^5+542x^4+972x^3+$\\
            &           &                                          & $+833x^2+276x+13=0$\\
 \hline
$181$ & ${\Bbb Q}(\sqrt{181}) $ & $435$ & $x^7+23x^6+210x^5+974x^4+2441x^3+$\\
            &             &                                         & $+3234x^2+2030x+435=0$\\
 \hline
$197$ & ${\Bbb Q}(\sqrt{197})$  & $7$ & $x^8+24x^7+228x^6+1095x^5+2834x^4+$\\
             &        &                                          &$+3942x^3+2795x^2+925x+112=0$\\
\hline
\end{tabular}

\bigskip
\centerline{Table 1.}
\end{table}

\bibliographystyle{amsplain}


\end{document}